\documentclass{article}
\usepackage{amssymb, amsmath}

\providecommand{\abs}[1]{\lvert#1\rvert}
\providecommand{\bigabs}[1]{\bigl\lvert#1\bigr\rvert}
\providecommand{\biggabs}[1]{\biggl\lvert#1\biggr\rvert}

\providecommand{\norm}[1]{\lVert#1\rVert}
\providecommand{\bignorm}[1]{\bigl\lVert#1\bigr\rVert}

\def\XXint#1#2#3{{\setbox0=\hbox{$#1{#2#3}{\int}$}
\vcenter{\hbox{$#2#3$}}\kern-.5\wd0}}

\providecommand{\supp}{\text{supp}}

\title{The topology of the spectrum for\\ Gelfand pairs on Lie groups}
\author{Fabio Ferrari Ruffino}
\date{}

\begin{document}

\newtheorem{Theorem}{Theorem}
\newtheorem{Lemma}[Theorem]{Lemma}
\newtheorem{Corollary}[Theorem]{Corollary}
\newtheorem{Remark}[Theorem]{Remark}
\newtheorem{Def}{Definition}

\maketitle

\begin{abstract}
Given a Gelfand pair of Lie groups, we identify the spectrum with a suitable subset of $\mathbb{C}^{n}$ and we prove the equivalence between Gelfand topology and euclidean topology.
\end{abstract}

\section{Introduction}

Let $(G,K)$ be a Gelfand pair with $G$ a connected Lie group and $K$ a compact subgroup. The Gelfand spectrum $\Sigma$ of $L^{1}(G)^{\natural}$, the commutative convolution algebra of bi-$K$-invariant integrable functions on $G$, is identified, as a set, with the set of bounded spherical functions. The Gelfand topology on $\Sigma$ is, by definition, the weak-$*$ topology, which coincides with the topology of uniform convergence on compact sets.

Since $G$ is a connected Lie group, the spherical functions on $G$ are characterized as the joint eigenfunctions of the algebra $\mathbb{D}(G/K)$ of differential operators on $G/K$ invariant by left $G$-translation. Being this algebra finitely generated, we idenfity $\Sigma$ with a subset of $\mathbb{C}^{s}$ assigning to each function the $s$-tuple of its eigenvalues with respect to a finite set of generators. Hence one can define on $\Sigma$ also the Euclidean topology induced from $\mathbb{C}^{s}$. In this article we prove that the two topologies coincide.

\section{Gelfand pairs and spherical functions}

We briefly recall the general theory of Gelfand pairs, that can be found in \cite{Faraut}.

Let $(G,\cdot)$ be a locally compact group, with a fixed left Haar measure $dx$. Let $K \leq G$ be a compact subgroup with normalized Haar measure $dk$.

\begin{Def} A function $f: G \rightarrow \mathbb{C}$ is said to be \emph{bi-invariant under $K$} if it is constant on double cosets of $K$, i.e., if:
	\[f(k_{1}xk_{2}) = f(x)\quad \forall k_{1},k_{2}\in K,\, \forall x\in G
\]
\end{Def}

Let $C_{c}(G)^{\natural}$ (resp. $L^{1}(G)^{\natural}$) be the set of continuous compactly-supported (resp. $L^{1}$) functions $f: G \rightarrow \mathbb{C}$ that are bi-invariant under $K$. It is easy to verify that it is a subalgebra of $C_{c}(G)$ (resp. of $L^{1}(G)$) with respect to the convolution in $G$. 

\begin{Def} $(G,K)$ is said to be a \emph{Gelfand pair} if $C_{c}(G)^{\natural}$ is a commutative algebra.
\end{Def}

One can easily prove that $C_{c}(G)^{\natural}$ is dense in $L^{1}(G)^{\natural}$, therefore $C_{c}(G)^{\natural}$ is a commutative algebra if and only if $L^{1}(G)^{\natural}$ is.

\paragraph{}Given a function $\varphi \in C(G)$ (not necessarily compactly-supported), we consider the linear functional:
\begin{equation*}
\begin{split}
\chi_{\varphi}:\;&C_{c}(G) \rightarrow \mathbb{C}\\
&\chi_{\varphi}(f) = \int_{G}f(x)\varphi(x^{-1})dx
\end{split}
\end{equation*}

\begin{Def} A function $\varphi \in C(G)$, $\varphi \neq 0$, is said to be \emph{spherical} if it is bi-invariant under $K$ and $\chi_{\varphi}$ is a character of $C_{c}(G)^{\natural}$, i.e.:
	\[\chi_{\varphi}(f*g) = \chi_{\varphi}(f)\cdot \chi_{\varphi}(g)\quad \forall f,g\in C_{c}(G)^{\natural}
\]
\end{Def}

One proves that $\varphi$ is spherical if and only if:
\begin{equation}\label{FormulaSpherical}
	\int_{K}\varphi(xky)dk = \varphi(x)\varphi(y) \quad \forall x,y\in G
\end{equation}
(see \cite{Faraut} prop. I.3 p. 319). In particular, this implies that $\varphi(1_{G}) = 1$.

\begin{Theorem} The dual space of $L^{1}(G)^{\natural}$ is $L^{\infty}(G)^{\natural}$. In fact, every continuous functional on $L^{1}(G)^{\natural}$ has the form:
	\[\chi_{\varphi}: f \rightarrow \int_{G}f(x)\varphi(x^{-1})dx
\]
with $\varphi \in L^{\infty}(G)^{\natural}$ unique and such that $\norm{\chi_{\varphi}} = \norm{\varphi}_{\infty}$.
\end{Theorem}

\paragraph{Proof:} If $\varphi \in L^{\infty}(G)^{\natural}$, $\chi_{\varphi}$ is a continuous functional on $L^{1}(G)$, hence on its closed subspace $L^{1}(G)^{\natural}$, and $\norm{\chi_{\varphi}} \leq \norm{\varphi}_{\infty}$.

\paragraph{}For the converse, let $\chi$ be a continuous functional on $L^{1}(G)^{\natural}$. By the Hahn-Banach theorem, we can extend $\chi$ to all of $L^{1}(G)$ without alterating its norm. So, being $L^{\infty}(G)$ the dual of $L^{1}(G)$, we have:
	\[\chi(f) = \int_{G}f(x)\psi(x^{-1})dx
\]
for some $\psi \in L^{\infty}(G)$, with $\norm{\chi} = \norm{\psi}_{\infty}$. Let $\varphi = \psi^{\natural}$ be the radialization of $\psi$, i.e.:
\begin{equation}\label{Radialization}
	\psi^{\natural}(x) = \iint_{K \times K}\psi(k_{1}xk_{2})dk_{1}dk_{2}
\end{equation}
It is easy to see that $\varphi\in L^{\infty}(G)^{\natural}$, and $\norm{\varphi}_{\infty} \leq \norm{\psi}_{\infty} = \norm{\chi}$. Moreover, $\psi$ and $\varphi$ induce the same functional on $L^{1}(G)^{\natural}$. In fact, if $f \in L^{1}(G)^{\natural}$, we have:
\begin{equation}\label{EquivalenceNatural}
\begin{split}
\int_{G}f(x)\varphi(x^{-1})dx &= \int_{G}f(x)\iint_{K\times K}\psi(k_{1}x^{-1}k_{2})dk_{1}dk_{2}dx \\
&= \iint_{K\times K}\int_{G}f(x)\psi(k_{1}x^{-1}k_{2})dx dk_{1} dk_{2} \\
&= \iint_{K\times K}\int_{G}f(k_{2}xk_{1})\psi(x^{-1})dx dk_{1} dk_{2} \\
&= \int_{G}f(x)\psi(x^{-1})dx
\end{split}
\end{equation}
Hence every continuous functional on $L^{1}(G)^{\natural}$ has the form $\chi_{\varphi}$ for $\varphi \in L^{\infty}(G)^{\natural}$, with $\norm{\varphi}_{\infty} \leq \norm{\chi_{\varphi}}$. Since we have also proved that $\norm{\chi_{\varphi}} \leq \norm{\varphi}_{\infty}$, we can conclude that $\norm{\chi_{\varphi}} = \norm{\varphi}_{\infty}$.
We now prove that $\varphi$ is unique: by linearity of $\chi_{\varphi}$ in $\varphi$, we have to prove that $\chi_{\varphi} = 0 \Rightarrow \varphi = 0$. But $\chi_{\varphi} = 0 \Leftrightarrow \norm{\chi_{\varphi}} = 0 \Leftrightarrow \norm{\varphi}_{\infty} = 0 \Leftrightarrow \varphi = 0$.\\
$\square$

\begin{Theorem} \emph{(See \cite{Faraut} Th. I.5 p. 320 or \cite{Helg1} Lemma 3.2 p. 408)} An element $\varphi$ of $L^{\infty}(G)^{\natural}$ defines a character of $L^{1}(G)^{\natural}$ if and only if $\varphi$ is a bounded spherical function.\\
$\square$
\end{Theorem}

\begin{Corollary}\label{BoundedSpherical} A bounded spherical function has $\infty$-norm equal to $1$.
\end{Corollary}

\paragraph{Proof:} If $\varphi \in L^{\infty}(G)^{\natural}$ is spherical, it determines a character $\chi_{\varphi}$ of $L^{1}(G)^{\natural}$, which is a commutative Banach algebra, with $\norm{\chi_{\varphi}} = \norm{\varphi}_{\infty}$. Hence, $\norm{\varphi}_{\infty} = 1$.\\
$\square$

\paragraph{}Let $\Sigma$ be the spectrum of $L^{1}(G)^{\natural}$, i.e., for the previous theorem, the set of bounded spherical functions. We define the \emph{Fourier spherical transform} (see \cite{Faraut} p. 333):
\begin{equation*}
\begin{split}
\hat{f}:\;&\Sigma\rightarrow \mathbb{C}\\
&\hat{f}(\varphi)=\chi_{\varphi}(f)=\int_{G}f(x)\varphi(x^{-1})dx
\end{split}
\end{equation*}
We can introduce on $\Sigma$ the \emph{Gelfand topology}, i.e., the weak-$*$ topology.

\begin{Theorem} The Gelfand topology on $\Sigma$ is equal to the topology of uniform convergence on compact sets (or locally uniform convergence).\\
$\square$
\end{Theorem}
(The proof is similar to the one given in \cite{Rudin1} p. 10-11.)

\section{The case of Lie groups}

If $G$ and $K$ are Lie groups, we can characterize Gelfand pairs and spherical functions by a differential point of view. Given a differential operator $D$ on a manifold $M$ (see \cite{Helg1} p. 239) and a diffeomorphism $\phi$ of $M$, we say that $D$ is $\phi$-invariant if $D(f\circ\phi) = Df \circ \phi\; \forall f\in C^{\infty}_{C}(M)$.

On a Lie group $G$ we have a special family of diffeomorphisms, the left translations by elements of $G$: $\phi_{g}(x) = gx$. Remembering that a Lie group always admits an \emph{analytic} structure compatible with the operations, we can construct a unique analytic structure also on the space of left cosets $G/K$ (with the quotient topology) such that the $G$-action on $G/K$:
\begin{equation*}
\begin{split}
L:\;&G\times G/K\rightarrow G/K\\
&L(x,gK)=xgK
\end{split}
\end{equation*}
is analytic (see \cite{Helg2} p. 113).

Let $C^{\infty}_{K}(G)$ be the set of functions in $C^{\infty}(G)$ such that $f(xk) = f(x)\; \forall k\in K, g\in G$. We have an isomorphism of algebras between $C^{\infty}(G/K)$ and $C^{\infty}_{K}(G)$ given by the projection $\pi$.

\paragraph{}We consider three algebras of differential operators (see \cite{Helg1} p. 274-287 and \cite{Helg2} p. 389-398):
\begin{table}[h]
	\centering
		\begin{tabular}{lcl}
			$\mathbb{D}(G)$ & $=$ & \{diff. op. on $G$ invariant by left $G$-translation\}\\
			$\quad\mathbb{D}_{K}(G)$ & $=$ & \{diff. op. in $\mathbb{D}(G)$ invariant also by right $K$-translation\}\\
			$\mathbb{D}(G/K)$ & $=$ & \{diff. op. on $G/K$ invariant by left $G$-translation\}
		\end{tabular}
\end{table}

We also consider the algebra:
	\[\mathbb{D}_{K}^{K}(G) = \mathbb{D}_{K}(G) / A,\quad A = \{D\in \mathbb{D}_{K}(G):\, Df = 0\; \forall f\in C^{\infty}_{K}(G)\}
\]
We can think of $\mathbb{D}_{K}^{K}(G)$ as the algebra of differential operators in $\mathbb{D}_{K}(G)$ acting only on $C^{\infty}_{K}(G)$: infact, if $D$ and $E$ coincide on $C^{\infty}_{K}(G)$, we have $D - E \in A$.

One can prove that $\mathbb{D}_{K}^{K}(G) \cong \mathbb{D}(G/K)$, with the isomorphism given by the projection $\pi$ (see \cite{Helg2} lemma 2.2 p. 390).

\begin{Theorem}\label{GelfandDiff} \emph{(See \cite{Helg1} p. 485 ex. 13)} Let $G$ be a \emph{connected} Lie group and let $K$ be a compact subgroup. Then, $(G,K)$ is a Gelfand pair if and only if $\;\mathbb{D}^{K}_{K}(G)$ is a \emph{commutative} algebra.\\
$\square$
\end{Theorem}

\begin{Theorem}\label{SphericalDiff} \emph{(See \cite{Helg1} prop. 2.2 p. 400)} Let $(G,K)$ be a Gelfand pair of Lie groups and $f\in C(G)$. Then, $f$ is spherical if and ony if:
\begin{itemize}
	\item $f\in C^{\infty}(G)^{\natural}$;
	\item $f(1_{G})=1$;
	\item $f$ is an eigenfunction of all the operators in $\mathbb{D}_{K}^{K}(G)$:
	\[Df = \lambda_{D}f \quad \forall D\in \mathbb{D}^{K}_{K}(G)
\]
\end{itemize}
$\square$
\end{Theorem}

\paragraph{Remark:} The proof of the theorem shows that a spherical function is necessarily \emph{analytic}.

\paragraph{}It can be proved that, being $K$ compact, $\mathbb{D}^{K}_{K}(G)$ is a \emph{finitely-generated} algebra (see \cite{Helg2} cor. 2.8 p. 395 and th. 5.6 p. 421). Let $D_{1}, ..., D_{s}$ be generators. Of course, $\varphi$ is an eigenfunction of all the operators in $\mathbb{D}^{K}_{K}(G)$ if and only if it is an eigenfunction of the generators. In this way, we can associate to each spherical function the $s$-uple of eigenvalues $(\lambda_{1}, ..., \lambda_{s})$ with respect to the generators. We can also prove that this association is injective, because the analyticity implies that two spherical functions having the same eigenvalues $(\lambda_{1}, ..., \lambda_{s})$ must coincide (see \cite{Helg1} cor. 2.3 p. 402).

\section{The topology of the spectrum}

In this way, we identify a spherical function, and in particular a bounded one, with a point in $\mathbb{C}^{s}$. So we identify the spectrum $\Sigma$ of $L^{1}(G)^{\natural}$ with a subset $A \subseteq \mathbb{C}^{s}$. Now, on $\Sigma$ we have the Gelfand topology, and on $A$ we have the induced euclidean topology. It is natural to ask if these two topologies coincide.

\begin{Lemma}\label{Sequences} \emph{(See \cite{Dugundji} p. 218)} Let $X$ be a topological space such that every point has a countable fundamental system of neighborhoods. Then a subset of $X$ is closed if and only if it is sequentially closed.

In particular, two such topologies coincide if and only if they induce the same notion of convergence on sequences.\\
$\square$
\end{Lemma}

\begin{Lemma}\label{Separable} $L^{1}(G)^{\natural}$ is separable.
\end{Lemma}

\paragraph{Proof:} We have only to prove that $L^{1}(G)$ is separable, because a subset of a separable space is separable (see \cite{Brezis} prop. III.22 p. 47). To see this, we choose a denumerable base of $G$ (which exists by definition of differential manifold) and we consider the subspace generated by the characteristic functions of these base-sets. Then we can argue as for $L^{1}(\mathbb{R}^{n})$ (see \cite{Brezis} Th. IV.13 p. 62).\\
$\square$

\begin{Remark}\label{LaplaceBeltrami} \emph{Being $K$ compact, one can construct on $G/K$ a riemannian metric invariant by the left action of $G$: hence, Laplace-Beltrami operator $\Delta$ with respect to this metric is invariant by left $G$-translations (see \cite{Helg2} Prop. 2.1 p. 387), i.e., $\Delta \in \mathbb{D}(G/K)$. This implies that, if $\varphi$ is a spherical function, $\pi: G \rightarrow G/K$ is the projection and $\varphi^{\pi} = \varphi \circ \pi^{-1}$, then $\varphi^{\pi}$ is an eigenfunction of $\Delta$, which is an elliptic operator.}
\end{Remark}

\begin{Theorem} The induced euclidean topology on $A$ and the Gelfand topology on $\Sigma$ coincide under the bijection $\varphi \in \Sigma \longleftrightarrow (\lambda_{1}, ..., \lambda_{s}) \in A$.
\end{Theorem}

\paragraph{Proof:} Of course $A$ is a metric space, so every point of $A$ has a denumerable fundamental system of neighborhoods. By corollary \ref{BoundedSpherical}, $\Sigma \subseteq B\Bigl(\bigl(L^{1}(G)^{\natural}\bigr)'\Bigr)$ (where $B$ is the unit ball). Being $L^{1}(G)^{\natural}$ separable for lemma \ref{Separable}, the weak-$*$ topology is metrizable on the unit ball (see \cite{Brezis} Th. III.25 p. 48), in particular the Gelfand topology on $\Sigma$ is metrizable. So, applying lemma \ref{Sequences}, we have to prove that the two topologies we are considering induce the same notion of convergence.

\paragraph{}Let $\{\varphi_{n}\}_{n\in \mathbb{N}}$ be a sequence of spherical functions, and let $\mathbb{D}^{K}_{K}(G) = \langle D_{1},$ $\ldots, D_{s} \rangle$. Let, $\forall i\in \{1,..,s\}$:
\begin{equation*}
\begin{split}
&D_{i}\varphi_{n} = \lambda_{i,n}\varphi_{n}\\
&D_{i}\varphi = \lambda_{i}\varphi
\end{split}
\end{equation*}
We have to prove that if $\varphi_{n} \rightarrow \varphi$ locally uniformly, then $\lambda_{i,n} \rightarrow \lambda_{i}\; \forall i\in \{1,..,s\}$. But $D_{i}\varphi_{n}(1_{G}) = \lambda_{i,n}\varphi_{n}(1_{G}) = \lambda_{i,n}$ and similarly $D_{i}\varphi(1_{G}) = \lambda_{i}$. So, being $D_{1}, ..., D_{s}$ generators, we have to prove that:
	\[\varphi_{n} \rightarrow \varphi \text{ loc. unif.} \Rightarrow D\varphi_{n}(1_{G}) \rightarrow D\varphi(1_{G}),\quad \forall D \in \mathbb{D}^{K}_{K}(G)
\]

If $f$ is a spherical function, it is continuous and non-zero by hypotesis, so it is easy to construct a function $\rho \in C^{\infty}_{c}(G)$ such that:
	\[\int_{G}f(x)\rho(x)dx \neq 0
\]
(We have to choose a point $x_{0}\in G$ such that $f(x_{0}) \neq 0$, choose by continuity a neighborhood $U(x_{0})$ such that $\Re f$ or $\Im f$ has constant sign on $U$, and construct $\rho \geq 0$ such that $\supp(\rho) \subseteq U$ and $\rho(x_{0}) = 1$). So we have, by the formula \eqref{FormulaSpherical}:
\begin{equation*}
\begin{split}
f(x)\int_{G}f&(y)\rho(y)\,dy = \int_{G}\rho(y)\biggl(\int_{K}f(xky)\,dk\biggr)\, dy\\
&=\int_{K}\int_{G}\rho(y)f(xky)\,dydk = \int_{K}\int_{G}\rho(k^{-1}x^{-1}y)f(y)\, dydk\\
&=\int_{G}\biggl(\int_{K}\rho(k^{-1}x^{-1}y)dk\biggr)f(y)\,dy
\end{split}
\end{equation*}
Concretely, the last integral in $dy$ is not extended to all of $G$: indeed, the domain of integration is the set of $y$ such that $\exists k\in K: k^{-1}x^{-1}y \in \supp(\rho)$, i.e., $x \cdot K \cdot \supp(\rho)$, which is compact because the product in $G$ is continuous.

So, if we restrict $x$ to an open neighborhood $V$ of $1_{G}$ with $\overline{V}$ \emph{compact}, we can assume that, for all such $x$, the domain of integration is $\overline{V} \cdot K \cdot \supp(\rho)$. We put:
\begin{equation*}
\begin{split}
&C = \overline{V} \cdot K \cdot \supp(\rho)\\
&A = \frac{1}{\int_{G}f(y)\rho(y)\,dy}\\
&\psi(x, y) = \int_{K}\rho(k^{-1}x^{-1}y)dk
\end{split}
\end{equation*}
$C$ is compact, $A \neq 0$ and, being $K$ compact, $\psi(x,y) \in C^{\infty}(G\times G)$. So, in particular, $\psi(\cdot, y) \in C^{\infty}(V)\, \forall y\in C$. We have, for $D \in \mathbb{D}(G)$:
\begin{equation*}
\begin{split}
	&f\vert_{V}(x) = A \int_{C}\psi(x,y)f(y)\, dy\\
	&Df\vert_{V}(x) = A \int_{C} \bigl[D^{(x)}\psi(x,y)\bigr] f(y)\,dy\\
	&Df\vert_{V}(1_{G}) = A \int_{C} \eta_{D}(y)f(y)\,dy
\end{split}
\end{equation*}
with $\eta_{D}(y) = \bigl(D^{(x)}\psi(x,y)\vert_{x=1_{G}}\bigr)$. But $\eta_{D}(y)$ is a continuous function, indeed $\psi \in C^{\infty}(G \times G)$, so $D^{(x)}\psi(x,y) \in C^{\infty}(G \times G)$, and, composing with the immersion $y \rightarrow (1_{G}, y)$ we still obtain a $C^{\infty}(G)$ function. So, the restriction of $\eta_{D}$ to $C$ is still continuous.

\paragraph{}So, applying the previous formula to $\varphi_{n}$ and $\varphi$, we obtain:
\begin{equation*}
\begin{split}
&D\varphi_{n}(1_{G}) = A_{n} \int_{C_{n}} \eta_{n,D}(y)\varphi_{n}(y)\,dy\\
&D\varphi(1_{G}) = A \int_{C} \eta_{D}(y)\varphi(y)\,dy
\end{split}
\end{equation*}
But, by construction, we can suppose $\eta_{n,D} = \eta_{D}$: in fact, we can begin the construction with $\rho_{n} = \rho$. For this, being $\varphi_{n}(1_{G}) = \varphi(1_{G}) = 1$, we choose a neighborhood $U(1_{G})$ with compact closure such that $\Re\varphi\vert_{U} \geq \delta > 0$. Then, being by hypotesis $\varphi_{n}\vert_{U} \rightarrow \varphi\vert_{U}$ uniformly, we can suppose that $\Re\varphi_{n}\vert_{U} > 0\; \forall n\in\mathbb{N}$. So we take $\rho_{n} = \rho$ such that $\rho(1_{G}) = 1$ and $\rho = 0$ outside $U$. From this we deduce that $\eta_{n,D} = \eta_{D}$ and $C_{n} = C$.

\paragraph{}If $\varphi_{n} \rightarrow \varphi$ uniformly on compact sets, in particular uniformly on $C$, being $\eta_{D}$ continuous and hence bounded on $C$, we have that $\eta_{D}\cdot\varphi_{n} \rightarrow \eta_{D}\cdot\varphi$ uniformly on $C$. So $\int_{C} \eta_{D}(y)\varphi_{n}(y)\,dy \rightarrow \int_{C} \eta_{D}(y)\varphi(y)\,dy$. Moreover, $A_{n} \rightarrow A$, in fact:
\begin{equation*}
\begin{split}
\biggabs{\int_{G}\varphi_{n}(x)&\rho(x)dx - \int_{G}\varphi(x)\rho(x)dx} \leq \int_{G}\bigabs{\varphi_{n}(x) - \varphi(x)}\rho(x)dx\\
&= \int_{\supp(\rho)}\bigabs{\varphi_{n}(x) - \varphi(x)}\rho(x)dx \leq K\int_{\supp(\rho)}\bigabs{\varphi_{n}(x) - \varphi(x)}dx \rightarrow 0
\end{split}
\end{equation*}
So $D\varphi_{n}(1_{G}) \rightarrow D\varphi(1_{G})$.

\paragraph{}For the converse, we know that $\Sigma \subseteq B\Bigl(\bigl(L^{1}(G)^{\natural}\bigr)'\Bigr)$, which is compact for the weak-$*$ topology by the Alaoglu-Banach theorem (see \cite{Brezis} Th. III.15 p. 42). Being the Gelfand topology metrizable on $B\Bigl(\bigl(L^{1}(G)^{\natural}\bigr)'\Bigr)$, compactness is equivalent to compactness by sequences (see \cite{Checcucci} prop. 4.4 p. 188). We indicate with $\overset{\bullet}\rightarrow$ the convergence with respect to the euclidean topology on $A$. So let us suppose that $\{\varphi_{n}\}_{n \in \mathbb{N}} \subseteq \Sigma$ is such that $\varphi_{n} \overset{\bullet}\rightarrow \varphi$. By compactness, we can extract a convergent subsequence (with respect to the Gelfand topology) $\varphi_{n_{k}} \rightarrow \tilde{\varphi}$, with $\tilde{\varphi} \in B\Bigl(\bigl(L^{1}(G)^{\natural}\bigr)'\Bigr)$. But necessarily $\tilde{\varphi} \in \Sigma \cup \{0\}$: indeed, $\chi_{\varphi_{n_{k}}}(f*g) \rightarrow \chi_{\tilde{\varphi}}(f*g)$ by definition on Gelfand topology, but $\chi_{\varphi_{n_{k}}}(f*g) = \chi_{\varphi_{n_{k}}}(f) \cdot \chi_{\varphi_{n_{k}}}(g) \rightarrow \chi_{\tilde{\varphi}}(f) \cdot \chi_{\tilde{\varphi}}(g)$.

\paragraph{}By remark \ref{LaplaceBeltrami}, the functions $\varphi_{n}^{\pi}$ are solutions of the equation:
	\[(\Delta - \lambda_{\Delta, n})\varphi_{n}^{\pi} = 0
\]
with $\Delta$ elliptic. Morover, $\lambda_{\Delta, n} \rightarrow \lambda_{\Delta}$ with $\lambda_{\Delta}$ defined by $\Delta \varphi^{\pi} = \lambda_{\Delta} \varphi^{\pi}$. Choosing a local chart $(U, \xi)$ in the origin of $G/K$, we have, by \cite{Gilbarg} Th. 8.32 p. 210, with $\Omega = U$, $\Omega' \subseteq \overline{\Omega'} \subseteq \Omega$, $f = g = 0$, $\alpha = 0$ and denoting by $\norm{\cdot}_{s}$ the Sobolev norm of order $s$:
	\[\bignorm{(\varphi_{n}^{\pi}\circ\xi^{-1}) \vert_{\xi(\Omega')}}_{1} \leq C\bigl(\bignorm{(\varphi_{n}^{\pi}\circ\xi^{-1}) \vert_{\xi(\Omega)}}_{0}\bigr) = C
\]
and one can easily verify that $C$ is independent by $n$ because $\lambda_{\Delta, n} \rightarrow \lambda_{\Delta}$, hence the sequence $\{\lambda_{n}\}_{n \in \mathbb{N}}$ is bounded. This implies that the functions $\varphi_{n}^{\pi} \vert_{\Omega'}$, and in particular the functions $\varphi_{n_{k}}^{\pi} \vert_{\Omega'}$, are equicontinuous and, by Arzela-Ascoli theorem (see \cite{Rudin} Th. 11.28 p. 245), there is a subsequence $\varphi_{n_{k_{h}}}^{\pi} \vert_{\Omega'} \rightarrow \psi^{\pi}$ locally uniformly on $G/K$. It is easy to deduce from this that $\varphi_{n_{k_{h}}} \vert_{\pi^{-1}(\Omega')} \rightarrow \psi$ locally uniformly on $G$. In particular, $\psi(1_{G}) = 1$ because $\varphi_{n_{k_{h}}}(1_{G}) = 1\; \forall k\in \mathbb{N}$.
We can choose $\Omega'' \subseteq \overline{\Omega''} \subseteq \pi^{-1}(\Omega')$ neighborhood of $1_{G}$ such that $\Re\psi\vert_{\Omega''} \geq \delta > 0$: in particular, $\varphi_{n_{k_{h}}}\vert_{\Omega''} \rightarrow \psi\vert_{\Omega''}$ uniformly. Then, if $\chi_{\Omega''}$ is the characteristic function of $\Omega''$, we consider the function:
	\[\xi(x) = (\chi_{\Omega''})^{\natural}(x^{-1})
\]
with $(\chi_{\Omega''})^{\natural}$ defined according to formula \eqref{Radialization} pag. \pageref{Radialization}. We have:
\begin{equation*}
\begin{split}
	&\hat{\xi}(\varphi_{n_{k_{h}}}) = \int_{G}\xi(x)\varphi_{n_{k_{h}}}(x^{-1})dx = \int_{G}(\chi_{\Omega''})^{\natural}(x^{-1})\varphi_{n_{k_{h}}}(x^{-1})dx\\
	&\qquad\quad\; = \int_{G}(\chi_{\Omega''})^{\natural}(x)\varphi_{n_{k_{h}}}(x)dx = \int_{G}\chi_{\Omega''}(x)\varphi_{n_{k_{h}}}^{\natural}(x)dx \\
	&\qquad\quad\; = \int_{G}\chi_{\Omega''}(x)\varphi_{n_{k_{h}}}(x)dx = \int_{\Omega''}\varphi_{n_{k_{h}}}(x)dx\\
&\Re\bigl[\hat{\xi}(\varphi_{n_{k_{h}}})\bigr] \rightarrow \int_{\Omega''}\Re\psi(x)dx \geq \delta \abs{\Omega''} > 0
\end{split}
\end{equation*}
Hence, by definition of Gelfand topology, it is not possible that $\varphi_{n_{k}} \rightarrow 0$, so that $\varphi_{n_{k}} \rightarrow \tilde{\varphi} \in \Sigma$.

\paragraph{}But, for the first part of the theorem, it must be $\varphi_{n_{k}} \overset{\bullet}\rightarrow \tilde{\varphi}$, so $\tilde{\varphi} = \varphi$. Hence, we have proved that for every sequence $\varphi_{n} \overset{\bullet}\rightarrow \varphi$, we can find a subsequence $\varphi_{n_{k}} \rightarrow \varphi$ uniformly on compact sets. Let us suppose that $\varphi_{n} \nrightarrow \varphi$: then, we can find a compact set $C \subseteq G$, $\varepsilon > 0$ and a subsequence $\varphi_{n_{k}}$ such that $\underset{x\in C}\sup \abs{\varphi_{n_{k}}(x) - \varphi(x)} > \varepsilon\; \forall k \in \mathbb{N}$. But, of course, $\varphi_{n_{k}} \overset{\bullet}\rightarrow \varphi$, so, applying the previous argument, we can find a sub-subsequence $\varphi_{n_{k_{j}}} \rightarrow \varphi$ uniformly on compact sets, in particular uniformly on C: a contradiction.\\
$\square$

\paragraph{}From the proof of the previous theorem one can also conclude that:

\begin{Corollary}\label{AClosed} $A$ is closed in $\mathbb{C}^{s}$.
\end{Corollary}

\paragraph{Proof:} Let $\{z_{n}\}_{n \in \mathbb{N}} = \{(\lambda_{1,n}, \ldots, \lambda_{s,n})\}_{n \in \mathbb{N}}$ be a sequence in $A$, with $z_{n} \rightarrow z = (\lambda_{1}, \ldots, \lambda_{s}) \in \mathbb{C}^{s}$. Let $\varphi_{n} \in \Sigma$ be the spherical function associated to $z_{n}$. We have that $\lambda_{\Delta,n} = P(\lambda_{1,n}, \ldots, \lambda_{s,n})$ with $P$ polynomial, hence $\lambda_{\Delta,n} \rightarrow \lambda_{\Delta} = P(\lambda_{1}, \ldots, \lambda_{n})$: the sequence $\{\lambda_{\Delta,n}\}$ is then bounded, hence, arguing as in the proof of the theorem, we can extract a subsequence $\varphi_{n_{k}} \rightarrow \tilde{\varphi} \in \Sigma$. But necessarily $\varphi_{n_{k}} \overset{\bullet}\rightarrow \tilde{\varphi}$, hence $(\lambda_{1}, \ldots, \lambda_{s})$ is the point of $A$ assiociated to $\tilde{\varphi}$, so $z \in A$.\\
$\square$

\begin{Corollary}\label{ConvDeriv} If $\varphi_{n} \rightarrow \varphi$ in $\Sigma$ then $D\varphi_{n} \rightarrow D\varphi$ uniformly on compact sets for every differential operator $D$.
\end{Corollary}

\paragraph{Proof:} For every $x_{0} \in G$, we have, for $V$ neighborhood of $x_{0}$ with $\overline{V}$ compact, $C$ compact and $A_{n} \neq 0$:
\begin{equation*}
\begin{split}
	&\varphi_{n}\vert_{V}(x) = A_{n} \int_{C}\psi(x,y)\varphi_{n}(y)\, dy\\
	&D\varphi_{n}\vert_{V}(x) = A_{n} \int_{C} \bigl[D^{(x)}\psi(x,y)\bigr] \varphi_{n}(y)\,dy\\
	&D\varphi_{n}\vert_{V}(x) = A_{n} \int_{C} \eta_{D}(x,y)\varphi_{n}(y)\,dy
\end{split}
\end{equation*}
with $\eta$ continuous. Similarly:
	\[D\varphi\vert_{V}(x) = A \int_{C} \eta_{D}(x,y)\varphi(y)\,dy
\]
By continuity, $\eta_{D}$ is bounded on $\overline{V} \times C$, so in particular on $V \times C$, so, being $\varphi_{n}\vert_{C} \rightarrow \varphi\vert_{C}$ uniformly, we have that $\int_{C} \eta_{D}(x,y)\varphi_{n}(y)\,dy \rightarrow \int_{C} \eta_{D}(x,y)\varphi_{n}(y)\,dy$ uniformly on $x$. Moreover, $A_{n} \rightarrow A$, so $D\varphi_{n}\vert_{V} \rightarrow D\varphi\vert_{V}$ uniformly.\\
$\square$

\section{Acknowledgements}

I am grateful to professor Fulvio Ricci for the helpfulness and the professionalism he has always shown since I began to work with him.

\end{document}